\documentclass[11pt, oneside]{article}   	
\usepackage{geometry}                		
\usepackage{graphicx}				
\usepackage{amssymb}
\usepackage{amsmath}
\title{How to prove this polynomial always has integer values at all integers}
\author{Wilberd van der Kallen}
\date{September 2015}							
\begin{document}
\maketitle
\begin{abstract}
The following problem was posed by user ``Kevin'' on Mathoverflow.
How to prove this polynomial always has integer values at all integers?

$P_m(x)=\sum_{i=0}^{m}\sum_{j=0}^{m}\binom{x+j}{ j}\binom{x-1}{ j}\binom{j}{ i}\binom{m}{ i}\binom{i}{ m-j}\frac{3}{(2i-1)(2j+1)(2m-2i-1)}.$

\noindent
We provide an answer.
\end{abstract}

\parindent=0pt\parskip.5\baselineskip

So
$$P_m(x)=\sum_{i=0}^{m}\sum_{j=0}^{m}\binom{x+j}{ j}\binom{x-1}{ j}\binom{j}{ i}\binom{m}{ i}\binom{i}{ m-j}\frac{3}{(2i-1)(2j+1)(2m-2i-1)}.$$

Our task is to show it takes integer values on integers.

As Kevin explains at

 [question 209140](http://mathoverflow.net/q/209140)

$P_m(x)$ is an even polynomial of degree $2m$ and he could show that $xP_m(x)$ always has integer values at all integers.

Folowing Wadim Zudilin we put 
$$B_k(x)=\binom{x+k}{2k}+\binom{-x+k}{2k}.$$  

For $k\geq0$ the $B_k$ are even polynomials of degree $2k$ that take integer values on integers. One has $B_k(k)=1$ for $k\geq1$, but $B_0(0)=2$. Further $B_k(i)=0$ for $|i|<k$. 
So the matrix $$(B_k(i))_{0\leq k\leq m}^{0\leq i\leq  m}$$ is triangular.

Every even polynomial $f(x)$ of degree $2j$ is clearly a linear combination of $B_0,\dots,B_j$ and the coefficients are determined by $f(0),\dots,f(j)$.
When $f(0)=0$ it is actually a linear combination of  $B_1,\dots,B_j$. 

 Rewrite $P_m(x)$ as $$P_m(x)=\sum_k d(m,k)B_k(x)$$ with $d(m,k)\in\mathbb{Q}$. As explained by Kevin, $P_m(k)$ vanishes if 
 $m>2|k|-2\geq0$ because all terms in the sum vanish. It can also be shown that  $P_m(0)=0$ for $m\geq2$, but that is more tricky. 
 Indeed we will show that $d(m,0)=0$ for $m\geq2$.
 
 Note that $P_m(x)$ visibly lies in the local ring $\mathbb{Z}_{(2)}$ for integer $x$.
So it suffices to show that $d(m,k)$ lies in $\mathbb{Z}_{(p)}$ for any  odd prime $p$. 
In fact we will find that the $d(m,k)$ are integers for $m\geq1$. And $d(0,0)=3/2$ lies in $\mathbb{Z}_{(p)}$ for our odd prime $p$.
 For $m$ not too large one may simply compute all $d(m,k)$.
The matrix $$(d(m,k))_{0\leq m\leq10}^{0\leq k\leq10}$$ looks like this
$$\left(
\begin{array}{ccccccccccc}
 \frac{3}{2} & 0 & 0 & 0 & 0 & 0 & 0 & 0 & 0 & 0 & 0 \\
 1 & -2 & 0 & 0 & 0 & 0 & 0 & 0 & 0 & 0 & 0 \\
 0 & 0 & 6 & 0 & 0 & 0 & 0 & 0 & 0 & 0 & 0 \\
 0 & 0 & 0 & 24 & 0 & 0 & 0 & 0 & 0 & 0 & 0 \\
 0 & 0 & 0 & 4 & 118 & 0 & 0 & 0 & 0 & 0 & 0 \\
 0 & 0 & 0 & 0 & 60 & 696 & 0 & 0 & 0 & 0 & 0 \\
 0 & 0 & 0 & 0 & 12 & 720 & 4824 & 0 & 0 & 0 & 0 \\
 0 & 0 & 0 & 0 & 0 & 336 & 8288 & 38240 & 0 & 0 & 0 \\
 0 & 0 & 0 & 0 & 0 & 60 & 6516 & 95928 & 336822 & 0 & 0 \\
 0 & 0 & 0 & 0 & 0 & 0 & 2520 & 109872 & 1131732 & 3215544 & 0 \\
 0 & 0 & 0 & 0 & 0 & 0 & 392 & 67904 & 1735320 & 13647840 & 32651544 \\
\end{array}
\right).$$
We will tacitly use it to deal with small values of $m$.

We will study the set $$V_p=\{ (m,k)\in \mathbb{Z}\times \mathbb{Z} \mid d(m,k)\in \mathbb{Z}_{(p)}\}.$$

  Using a method of Zeilberger we will prove relations between the $d(m,k)$ that were first discovered experimentally.
 One relation allows us to rewrite  $m(m-1)(1 + 2 m) d(m, k)$ in such a manner that we can use
 the method of Floors described in 
 
 [question 26336](http://mathoverflow.net/q/26336).
 
 With that method we show that $m(m-1)(1 + 2 m) d[m, k]$ is an integer multiple of $3m(m-1)$.
 Together with the relations this will allow us to show that $V_p$
 fills all of $ \mathbb{Z}\times \mathbb{Z}$ for odd primes $p$.

Our variables 
$i,j,k,m,n,q$ will 
take integer values only.

As in the A=B book \cite{A=B} we use the convention that $\binom{x}{j}$ is a polynomial in $x$ for fixed $j$. And it is the zero polynomial if $j<0$. So $\binom i j$ is defined for all integers $i$, $j$. It also vanishes  if $j>i\geq0$. Of course $\binom{i}{j}$  agrees with the usual binomial coefficient if $0\leq j\leq i$.

By inspecting the values at $x=0,\dots,j,$ we see that $$(-1)^j\binom{x+j}{j}\binom{x-1}{j}-(-1)^{j-1}\binom{x+j-1}{j-1}\binom{x-1}{j-1}$$ equals $(-1)^j\binom{2j}{j}B_j(x)/2$ for $j\geq0$.
Taking the telescoping sum over $j$ gives $$(-1)^j\binom{x+j}{j}\binom{x-1}{j}=\sum_{k=0}^j(-1)^k\binom{2j}{j}B_k(x)/2$$ for $j\geq0$. (Valid for all $j$, actually).

This allows us to conclude that
$$d(m,k)=\sum_{i=0}^m\sum_{j=k}^m\frac{ 3(-1)^{k+j}\binom{2k  }{ k }\binom{j  }{  i}\binom{ m }{ i }\binom{i  }{ m-j }}{2(2i-1)(2j+1)(2m-2i-1)}.$$

In particular $d(m,k)=0$ for $m<0$ and for $k>m$. We will see that $m(m-1)d(m,k)$ also vanishes for
$2k-2<m$. 

Let us use the notation $[$statement$]=
\begin{cases}1,&\text{if statement is true;}\\ 0,&\text{otherwise.}\end{cases}$

Then 
\begin{align}
d(m,k)=&\sum_{i,j}[j\geq k\geq0] \text{term}(m,k,i,j), \tag{$\Sigma ij$}
\end{align}

where $$\text{term}(m,k,i,j)=[m\geq0]\frac{3(-1)^{k+j}\binom{2k  }{ k }\binom{j  }{  i}\binom{ m }{ i }\binom{i  }{ m-j }}{2(2i-1)(2j+1)(2m-2i-1)}.$$

Put

\begin{align*}
\text{rel1}&(m,k)
 = \\
   -& 32 (3-2 k)^2 (-k+m+1)
   (-k+m+2) d(m,k-2)\\
 +&4 (-k+m+1) \left(2 k m^2-2 (k-1) (8 k-9) m+(2 k-3) (8
   (k-2) k+9)\right) d(m,k-1)\\
   + & k (-2 k+m+2) (-2 k+m+3) (-2 k+2
   m+1) d(m,k),
\end{align*}   
   
\begin{align*}
\text{rel2}&(m,k)
 = \\&-4 \left((m-1)^2-1\right) d(m-1,k-1)\\&-4 (2 (k-1)+m+1)
   (-k+m+1) d(m,k-1)\\&+k (2 k-m-2) d(m,k)
 \end{align*}   

\noindent{\bf Key results}

$\bullet$ $\text{rel1}(m,k)$  vanishes.

$\bullet$ $m(m-1)d(m,k)$ vanishes for $2k-2<m$. 

$\bullet$ $\text{rel2}(m,k)$ vanishes.

$\bullet$ $m(m-1)( 2 m+1) d(m, k)$ is an integer multiple of $3m(m-1)$.

Before proving the Key results, let us draw conclusions from them. Let $m\geq2$. As $d(m,0)=0$, we have $P_m(0)=0$ and the $d(m,k)$ are 
determined by $P_m(1)\dots,P_m(m)$. 
Now the integral matrix $$(B_k(i))_{1\leq k\leq m}^{1\leq i\leq  m}$$ is triangular with ones on the diagonal.
We conclude that $d(m,k)\in \mathbb{Z}_{(2)}$ for $m\geq2$.

Let $p$ be a prime, $p\geq5$, and let $m\geq 2$.
If $p$ does not divide $2m+1$, then $d(m,k)\in \mathbb{Z}_{(p)}$ because $m(m-1)( 2 m+1) d(m, k)\in 3m(m-1)\mathbb{Z}_{(p)}$.
Now assume $p$ divides $2m+1$. Then it does not divide $2m+3$, so then $d(m+1,j)\in \mathbb{Z}_{(p)}$ for all $j$.
Also, $p$ does not divide $(m-1)(m+1)$, so it follows from $\text{rel2}(m+1,k+1)=0$ that $d(m,k)\in \mathbb{Z}_{(p)}$.
We have shown that  $d(m,k)\in \mathbb{Z}_{(p)}$ if  $p$ is prime, $p\geq5$, $m\geq2$. 

Remains $p=3$. Let $m\geq2$ again.

If $3$ does not divide $2m+1$, then $d(m,k)\in 3\mathbb{Z}_{(3)}$ because $m(m-1)( 2 m+1) d(m, k)\in 3m(m-1)\mathbb{Z}_{(3)}$.

If $m\equiv1$ mod 9, or $m\equiv7$ mod 9, then $(2m+1)/3$ is prime to $3$ and $d(m,k)\in \mathbb{Z}_{(3)}$ because $m(m-1)(( 2 m+1)/3 )d(m, k)\in m(m-1)\mathbb{Z}_{(3)}$.

If $m\equiv 4$ mod 9, then $(m-1)(m+1)/3$ is prime to 3 and $d(m,k)\in \mathbb{Z}_{(3)}$ because $\text{rel2}(m+1,k+1)=0$ shows 
$((m-1)(m+1)/3)d(m,k)$ is an integer linear combination of the integers $d(m+1,j)/3$.

We conclude that $d(m,k)\in \mathbb{Z}_{(3)}$ for $m\geq2$.
So the $d(m,k)$ are integers for $m\geq2$ and $P_m$ takes integer values on integers for $m\geq2$.
Recall that $P_0$, $P_1$ also take integer values. $\text{\bf Done}$.

{\bf So we still have to prove the Key results}.

First a technical issue.
If $x>0$ then $\binom{x}{j}=\frac{\Gamma(1+x)}{\Gamma(1+j)\Gamma(1+x-j)}$ and the bimeromorphic function $$f(x,y)=\frac{\Gamma(1+x)}{\Gamma(1+y)\Gamma(1+x-y)}$$ is continuous at $(x,j)$. However, if $i<0$ then $f$ has an indeterminate value at $(i,j)$. For example, $\binom{i}{i}$ equals $1$ if $i\geq0$, but it vanishes for $i<0$. At $(-1,-1)$ both $0$ and $1$ are values of $f$. Indeed Mathematica can be steered to give either answer.

$\mathtt {Binomial[i,j]~ /.~ i->-1~/.~j->-1}$ gives 1 
and 

$\mathtt {Binomial[i,j]~ /. ~j->-1~/.~i->-1}$ gives 0.

And  $\mathtt{FullSimplify[Binomial[i, i] == Binomial[i - 1, i - 1]]}$ yields $\mathtt {True}$.
This answer is correct, but it tells only that for generic complex numbers $i$ the
identity holds.

Thus we need to make case distinctions when using identities between multimeromorphic functions, explicitly or implicitly, to prove identities involving
the $\binom{i}{j}$.

We start proving that $\text{rel1}(m,k)$ vanishes.

As $[j\geq k+1]\big(2 (2 k+1)\text{term}(m,k,i,j)+(k+1)\text{term}(m,k+1,i,j)\big)=0$, we get from $(\Sigma ij)$ that

\begin{align}
2 (2 k+1) d(m,k)+(k+1) d(m,k+1)=&\sum_i \text{iterm}(m,k,i)\tag{$\Sigma i$}
\end{align}
where $$\text{iterm}(m,k,i)=2 (2 k+1)\text{term}(m,k,i,k).$$

Now we use the

{\tt    Fast Zeilberger Package version 3.61\\
    written by Peter Paule, Markus Schorn, and Axel Riese\\
    Copyright 1995-2015, Research Institute for Symbolic Computation (RISC),\\
    Johannes Kepler University, Linz, Austria.}

It suggests to put

$$g(m,k,i)=\frac{3\times 2^{2 k+3} m (-2 i+m+1) \Gamma \left(k+\frac{3}{2}\right) \binom{k+1}{i-1}
   \binom{m-1}{k+1} \binom{k+1}{m-i}}{\Gamma \left(\frac{1}{2}\right) \Gamma(k+2)}$$

and show that
\begin{align*}
-32 (1 + 2 k) &(3 + 2 k) (k - m) (1 + k - m) \text{iterm}(m, k, i) \\
   -4 (1 + k - m) & (57 + 110 k + 72 k^2 + 16 k^3 - 34 m - 46 k m - 
   16 k^2 m + 4 m^2 + 2 k m^2)\\&\times \text{iterm}(m, k+1, i)\\
-(2 + k) (5 +  &2 k - 2 m) (3 + 2 k - m) (4 + 2 k - m)\text{iterm}(m, k+2, i)
 \\
  & -g(m, k, i + 1) +g(m, k, i)=0
\end{align*}   

for $m\geq0$.
So we do that and then sum over $i$, using $(\Sigma i)$. The $g$ terms drop out by telescoping and we get a relation 
\begin{align*}
-32 (1 + 2 k) &(3 + 2 k) (k - m) (1 + k - m) (2 (2 k+1) d(m,k)+(k+1) d(m,k+1) )\\
   -4 (1 + k - m) & (57 + 110 k + 72 k^2 + 16 k^3 - 34 m - 46 k m - 
   16 k^2 m + 4 m^2 + 2 k m^2)\\&\times  (2 (2 k+3) d(m,k+1)+(k+2) d(m,k+2) )\\
-(2 + k) (5 +  &2 k - 2 m) (3 + 2 k - m) (4 + 2 k - m) \\\times(2 (2 k+5) &d(m,k+2)+(k+3) d(m,k+3) )
 \\
  & =0
\end{align*}   
valid for all $m$, as it is obvious for $m<0$.
We may rewrite it as a recursion for $\text{ rel1}$:
$$2 (3 + 2 k) \text{rel1}(m, k + 2) + (2 + k)\text{ rel1}(m, k + 3)=0.$$
As $d(m,k)$ vanishes for $k>m$, it follows from the recursion that $\text{rel1}(m, k )$ vanishes for all $k$.

{\bf So we have established the vanishing of } $\text{rel1}(m, k )$.

Put
$$\text{pterm}(m,x,i,j)=\binom{x+j}{ j}\binom{x-1}{ j}\binom{j}{ i}\binom{m}{ i}\binom{i}{ m-j}\frac{3}{(2i-1)(2j+1)(2m-2i-1)},$$
so that 
$$P_m(x)=\sum_{i,j}\text{pterm}(m,x,i,j).$$
If $k\geq1$ and $\text{pterm}[m,k,i,j]$ is nonzero, then 
$ k-1\geq j$ and $m\geq j \geq i \geq m-j $. We see that $$P_m(k)=0 \text{ if }0\leq 2k-2<m,$$ because all the $\text{pterm}(m,k,i,j)$ vanish.
In particular we get
$$0=P_m(1)=\sum_kd(m,k)B_k(1)=2d(m,0)+d(m,1),$$ and $$0=P_m(2)=\sum_kd(m,k)B_k(2)=2d(m,0)+4d(m,1)+d(m,2)$$ for $ m\geq 3$.
So then $d(m,1)=-2d(m,0)$ and $d(m,2)=6d(m,0)$. Substitute this into $\text{rel1}(m, 2 )=0$ and you find 
$$4 m (m - 1) (-2 m - 1) d(m, 0)=0.$$
This means that $d(m, 0)=0$ for $m\geq3$. As $d(2,0)$ also vanishes, we now know that $m(m-1)P_m(k)$ vanishes
if $m>2|k|-2$.
As the matrix $$(B_k(i))_{0\leq k\leq m}^{0\leq i\leq  m}$$ is triangular, we now  conclude that
\begin{align}
m(m-1)d(m,k) \text{ vanishes for }&m>2k-2. \tag{SSE}
\end{align}

{\bf So we have established the vanishing of }   
$m(m-1)d(m,k)$ {\bf{ for }} $m>2k-2$.

Before turning to  $\text{rel2}(m, k )$ we compute $d(2k-2,k)$ and $d(2k-3,k)$ for $k\geq3$. These are the values that help to compute all $d(m,k)$ 
recursively with the recursion given by $\text{rel1}(m, k )$=0.
As $d(2k-2,j)$ vanishes for $j<k$, one has $$d(2k-2,k)=P_{2k-2}(k)=\text{pterm}(2k-2,k,k-1,k-1)$$
and similarly
$$d(2k-3,k)=P_{2k-3}(k)=\text{pterm}(2k-3,k,k-1,k-2)+\text{pterm}(2k-3,k,k-1,k-1).$$
So we know $d(m,k)$ for $m\geq 2k-3\geq3$. By (SSE) we also know  $d(m,k)$  for $k\leq1$ and any $m$.
Using these values we get $\text{rel2}(m,k)=0$ by inspection for $m\geq 2k-3$ or $k\leq1$.
Notice that $(-7+2 k)\text{rel2}(2k-4,k)-\text{rel1}(2k-4,k)$ is a combination of the known terms $ d(-5 + 2 k, -1 + k)$, 
$    d(-4 + 2 k, -2 + k)$,  $  d(-4 + 
      2 k, -1 + k)$. It also vanishes by inspection, so we now have that $\text{rel2}(m,k)=0$ for $m\geq 2k-4$ or $k\leq1$.

By substituting the definitions and expanding we check that
\begin{align*}
&(-1 + k) (-1 + 2 k - 2 m) (-3 + 2 k - m)(4 - 2 k + m) \text{rel2}(m, k)
    \\& + 
 4 (-1 + (-1 + m)^2) \text{rel1}(m - 1, k - 1)
 \\& - 
 32 (5 - 2 k)^2 (-2 + k - m) (-1 + k - m) \text{rel2}(m, k - 2)
 \\& + 
 4 (1 - k + m)\\&\times (-99 + 16 k^3 - 2 m (32 + m) - 8 k^2 (11 + 2 m) + 
    2 k (81 + m (31 + m))) \text{rel2}(m, k - 1)
   \\& - (-1 + k) (-4 + 2 k - m)\text{ rel1}(m, k)\\&  - 
 4 (-1 + k - m) (-5 + 2 k + m) \text{rel1}(m, k - 1)
 \\&=0
 \end{align*}   

As $\text{rel1}$ vanishes, this leads to the following recursion for $\text{rel2}$.
\begin{align*}
&(-1 + k) (-1 + 2 k - 2 m) (-3 + 2 k - m)(4 - 2 k + m) \text{rel2}(m, k) \\& - 
 32 (5 - 2 k)^2 (-2 + k - m) (-1 + k - m) \text{rel2}(m, k - 2)
 \\& + 
 4 (1 - k + m)\\&\times (-99 + 16 k^3 - 2 m (32 + m) - 8 k^2 (11 + 2 m) + 
    2 k (81 + m (31 + m))) \text{rel2}(m, k - 1)
   \\&=0
 \end{align*}  

As $\text{rel2}(m,k)=0$ for $ 2k-4\leq m$ or $k\leq1$,  the recursion shows by induction on $k$  
that $\text{rel2}(m,k)=0$ for all $m$, $k$.
    
    {\bf So we have also established the vanishing of } $\text{rel2}(m,k)$ and it is time to
    show the Key result that $m(m-1)( 2 m+1) d(m, k)$ is an integer multiple of $3m(m-1)$. 
This is obvious for $m<2$, so we further assume $m\geq2$. Then we know that $d(m,0)=0$ and we have seen this implies
$d(m,k)\in \mathbb{Z}_{(2)}$. So it suffices to
show that $m(m-1)( 2 m+1) d(m, k)\in 3m(m-1)\mathbb{Z}[1/2]$.

Using relation $(\Sigma i)$ we may rewrite $\text{rel1}(m,k)=0$ as
\begin{align*}&2 (m-1 ) m ( 2 m+1) d(m, k-1 )\\&+
(2 - 2 k + m) (3 - 2 k + m) (1 - 2 k + 2 m) \sum_i \text{iterm}(m, k - 1,i)\\&+ 
 16 (3 - 2 k) (-2 + k - m) (-1 + k - m) \sum_i \text{iterm}(m, k - 2,i) \\&=0
 \end{align*}  
 
We claim that 
\begin{align*}&(2 - 2 k + m) (3 - 2 k + m) (1 - 2 k + 2 m)  \text{iterm}(m, k - 1,i)\\&+ 
 16 (3 - 2 k) (-2 + k - m) (-1 + k - m)  \text{iterm}(m, k - 2,i) 
 \end{align*} 
 lies in $3m(m-1)\mathbb{Z}[1/2]$. 
  
That will prove that the $ (m-1 ) m ( 2 m+1) d(m, k-1 )$ are  integer multiples of  $3m(m-1)$.

Put 
$$\text{frac1}(m,k,i)=\frac{3 (m-1) m \binom{2 (k-1)}{k-1} (-2 k+2
   m+1) \binom{k-1}{i} \binom{m}{i} \binom{i}{-k+m+1}}{(2
   i-1) (2 m-2 i-1)}$$
   and 
   $$\text{frac2}(m,k,i)=6  (k-m-1)\binom{2 (k-1)}{k-1}  \binom{k-1}{i} \binom{m}{i} \binom{i}{-k+m+1}.$$
   
   Then $\text{frac1}(m,k,i)+\text{frac2}(m,k,i)$ equals 
   \begin{align*}&(2 - 2 k + m) (3 - 2 k + m) (1 - 2 k + 2 m)  \text{iterm}(m, k - 1,i)\\&+ 
 16 (3 - 2 k) (-2 + k - m) (-1 + k - m)  \text{iterm}(m, k - 2,i) ,
 \end{align*} 
 
 so it suffices to show that $\text{frac1}(m,k,i)/(6m(m-1))$ and $\text{frac2}(m,k,i)/(6m(m-1))$, which make sense for $m\geq2$,
  lie in $\mathbb{Z}[1/2]$ for $m\geq 2$.
 Recall that the Catalan numbers
 $$C(i)=\frac{\binom{2 i}{i}}{i+1}$$ are integers. See
 
 [A000108](https://oeis.org/A000108)
 
We now look at $\text{frac1}(m,k,i)/(6m(m-1))$.

If $\text{frac1}(m,k,i)$ is nonzero then $ m\geq k-1\geq i\geq m+1-k\geq 0$. 
We distinguish two cases: $ m=k-1\geq i\geq 0$ and $ m> k-1\geq i\geq m+1-k\geq 0$. 

First let $ m=k-1\geq i\geq 0$. If $i=k-1$, then 
$$\text{frac1}(m,k,i)/(6m(m-1))=\text{frac1}(k-1,k,k-1)/(6(k-1)(k-2))=C(k-2).$$

Similarly $\text{frac1}(k-1,k,0)/(6(k-1)(k-2))=C(k-2).$

So we may assume $0<i<m=k-1$.
Then

$\text{frac1}(m,k,i)/(6m(m-1))=\text{frac1}(m,m+1,i)/(6m(m-1))$ equals

$$
\frac{-(2 i-2)! (2 m)! (-2 i+2 m-2)!}{2 (i!)^2 (2 i-1)!
   ((m-i)!)^2 (-2 i+2 m-1)!}$$ and we must show it takes values in $\mathbb{Z}[1/2]$.

This is the kind of expression to which one may apply the method of Floors explained in 
 
  [question 26336](http://mathoverflow.net/q/26336).

It is based on
$$\operatorname{ord}_p n!
=\biggl\lfloor\frac{n}{p}\biggr\rfloor+\biggl\lfloor\frac{n}{p^2}\biggr\rfloor
+\biggl\lfloor\frac{n}{p^3}\biggr\rfloor+\dots$$

According to the method it suffices to check that $\text{test}(m,i,2n+1)\geq0$ for $n\geq1$, where
\begin{align*}
\text{test}(m,i,q)=&\\-2& \left\lfloor \frac{m-i}{q}\right\rfloor +\left\lfloor
   \frac{-2 i+2 m-2}{q}\right\rfloor -\left\lfloor \frac{-2
   i+2 m-1}{q}\right\rfloor\\ -2& \left\lfloor
   \frac{i}{q}\right\rfloor +\left\lfloor \frac{2
   i-2}{q}\right\rfloor -\left\lfloor \frac{2
   i-1}{q}\right\rfloor +\left\lfloor \frac{2
   m}{q}\right\rfloor
.\end{align*}

This is a tedious puzzle. For fixed $q$ the function $\text{test}(m,i,q)$ is periodic of period $q$ in both variables $i$ and $m$.
So for fixed $q$ one may simply compute all values. We do it for $3\leq q=2n+1<17$. The results are nonnegative.
But if $q$ is large we need to be more efficient. If both $q=2n+1$ and $m$ are
fixed, then $\text{test}(m,i,q)$ can only change value where at least one of the Floors jumps as a function of $i$. So it suffices to sample around the jumping
points (modulo $q$). We know where they are. More specifically, we only need to consider the 15 cases where one of $i-1$, $i$, $i+1$ lies in 
$\{0, 1, -1 + m, m,  -2 + m - n, -1 + m - n, 1 + n\}$.
So we can eliminate $i$ at the expense of having 15 cases. Similarly we can eliminate $m$ for each of those cases, ending up with
153 test functions that depend on $n$ only. Each test function is a linear combination of seven Floors.
Each of the Floors stabilises after $n$ has reached an easily computable bound.
For instance $\left\lfloor -\frac{8}{2 n+1}\right\rfloor$ is constant for $n\geq4$.
In fact the bound 5 suffices for all $7\times 153$ Floors. Compute the 153 stable values. They are nonnegative. 
 This solves the puzzle; the check for $3\leq q=2n+1<17$ was overkill.

 So we now turn to the case $ m> k-1\geq i\geq m+1-k\geq 0$.
 Then 
 \begin{align*}&\text{frac1}(m,k,i)/(6m(m-1)C(i-1))=\\&\frac{i! (2 k-2)! m! (-2 i+2 m-2)! (-2 k+2 m+1)!}{(2 i)! (k-1)! (-i+k-1)! (m-i)! (-2 i+2 m-1)!
   (-k+m+1)! (2 m-2 k)! (i+k-m-1)!}
 \end{align*} 
 
 We use  the method of Floors again to show that  $\text{frac1}(m,k,i)/(6m(m-1)C(i-1))\in \mathbb{Z}[1/2]$.
This time we eliminate $k$, $m$, $i$ in that order and take $n\geq6$ as bound where all  $13\times 3508$ Floors are stable.

{\bf So we have shown that} $\text{frac1}(m,k,i)/(6m(m-1))$ 
  lies in $\mathbb{Z}[1/2]$ for $m\geq 2$.
Remains showing that $\text{frac2}(m,k,i)/(6m(m-1))$ 
  lies in $\mathbb{Z}[1/2]$ for $m\geq 2$.

If $\text{frac2}(m,k,i)$ is nonzero then $m> k-1\geq i\geq m+1-k>0$ and  
$\text{frac2}(m,k,i)/(6m(m-1))$ 
equals 
$$\frac{-(2 k-2)! (m-2)!}{i! (k-1)! (-i+k-1)! (m-i)! (m-k)!
   (i+k-m-1)!}.$$
This can be treated like the previous case. We eliminate $k$, $m$, $i$ in that order and take $n\geq6$ as bound where all  $8\times 1278$ Floors are stable.

{\bf Done}

 \end{document}